\documentclass[10pt]{article}

\usepackage{amsfonts}
\usepackage{amsmath}
\usepackage{graphicx}

\begin{document}

\title{Admissible surfaces in progressive addition lenses}

\author{Sergio Barbero\footnote{Instituto de \'Optica (CSIC), Serrano 121, Madrid, 28006, Spain}\quad and\quad
Mar\'ia del Mar Gonz\'alez\footnote{Departamento de Matem\'aticas, Universidad Aut\'onoma de Madrid, Madrid, 28049, Spain}}

\maketitle




\begin{abstract}
Progressive addition lenses contain a surface of spatially-varying curvature, which provides variable optical power for different viewing areas over the lens. We derive complete compatibility equations that provide the exact magnitude of cylinder along lines of curvatures on any arbitrary PAL smooth surface. These equations reveal that, contrary to current knowledge,  cylinder, and its derivative, does not only depend on principal curvature and its derivatives along the principal line but also on the geodesic curvature and its derivatives along the line orthogonal to the principal line.
We quantify the relevance of the geodesic curvature through numerical computations. We also derive an extended and exact Minkwitz theorem only restricted to be applied along lines of curvatures, but excluding umbilical points.
\end{abstract}


\section{Introduction}

The use of progressive addition lenses (from now on PALs) is a widespread ophthalmic solution when different refractive correction is required for different object locations: far, near and intermediate.
PALs contain a surface of spatially-varying curvature, which provides the required variable optical power for the different viewing areas over the lens. However, the major drawback of PAL technology is the impossibility, intrinsic to any smooth surface, of providing spatially-varying mean curvature without introducing undesired astigmatism (differences between principal curvatures), except along umbilical lines; a fact firstly pointed out in PAL technology by Minkwitz~\cite{Minkwitz}.
This limitation imposes to PAL design  a trade-off between the required
mean curvature and the undesired astigmatism.
A way to tackle such trade-off is to minimize a functional where the function space $S$
is defined as all possible surfaces ($u \in S$) describing the variable curvature surface, and whose image
is the metric related to the aforementioned trade-off~\cite{Wang}:
\begin{equation}\label{functional}
J[u] = \int_{\Omega} \left\{ \alpha C(x,y)^2 + \beta(H(x,y)-H_{t}(x,y))^2\right\}dxdy
\end{equation}
where, $\Omega$ defines the domain of $u$; $C$, $H$ are cylinder and mean curvature
functions, respectively; $H_t$ is the prescribed mean curvature function; and $\alpha$ and $\beta$ are arbitrarily weight functions setting \textit{a priori} trade-off between the desired mean curvature and the undesired cylinder. Minimizing the functional in \eqref{functional} opposes serious problems unless good candidates for weight functions are selected~\cite{Alonso}.
However, \eqref{functional} could potentially benefit if constraints
limit the space $S$ in which we search for extremals of $J[u]$~\cite{Brunt}.
There are two type of constraints that could be imposed to \eqref{functional}. The first type is related to constraints in the manufacturing process; for instance, those imposed by the cutting tool radius on the maximum attainable curvatures~\cite{Alonso}. The second type is related to mathematical compatibility conditions between $C$ and $H$ for a surface to exist. In this letter, we concentrate in the second type.
Even without applying them to optimize \eqref{functional}, such relations are highly
relevant for a better understanding on how far PAL design can go in achieving its ideal target, i.e. $C(x,y)=0$ and $H(x,y)=H_{t}(x,y)$ for all points in the domain $\Omega$.

In PAL design, it is usually defined the so-called principal line as the curve embedded in the surface where it is prescribed the change of, at least, one principal curvature.
Among the surface compatibility conditions, to date only the so-called Minkwitz theorem~\cite{Minkwitz} and a generalization to it~\cite{Esser, Blendowske} have been found. Minkwitz theorem relates cylinder derivatives along $s_v$ with principal curvature derivatives along $s_u$, where $s_u$ and $s_v$ define arc-length along the principal and its orthogonal line, respectively: $\frac{d C}{d s_{v}}  = 2 \frac{d k_{u}}{d s_{u}}$.
However, Minkwitz theorem is only valid when the principal line is umbilic. A version of Minkwitz theorem, given directly the cylinder was obtained by Alonso et al~\cite{Alonso}: $C  = 2 \frac{d k_{u}}{d s_{u}} s_{v}$.

Trying to overcome this strong restriction of Minkwitz theorem,  Esser et al.~\cite{Esser} found new equations for non-umbilical lines, still under the restriction of being only applicable to principal lines of a symmetrical surface (hence, being planar curves). An alternative derivation to that of Esser et al. was provided by Blendowske~\cite{Blendowske}, avoiding the symmetry restriction. However, both  Minkwitz's theorem and the aforementioned generalizations are based on series expansion of either curvatures~\cite{Minkwitz, Esser} or the surface itself~\cite{Blendowske}, with a posterior truncation. As will be seen later, this approximation leads to a misleading idea: namely, the control of cylinder derivative depends locally \textit{only} on first derivatives of mean curvature (Minkwitz) and, additionally, on cylinder~\cite{Esser, Blendowske}  along exclusively the principal line. Loss of dependency on curvature properties along the orthogonal line to the principal one is a side effect of the Taylor truncation.

 In this work, we derive all possible exact equations relating cylinder and principal curvatures and/or with its derivatives along surface curves for any arbitrary surface.
The procedure for it is based on  surface existence theory, which establishes that Gauss-Codazzi-Mainardi equations (sometimes called compatibility equations) are necessary and sufficient surface differential conditions~\cite{Stoker}. Moreover, there is no way to derive other expressions relating first and second fundamental forms (p.239~\cite{Carmo}).
Any extra equation derived from these ones should be formally equivalent. This is indeed the case of Minkwitz theorem, as we will show, and was already hinted by Blendowske~\cite{Blendowske}.

\section{Theoretical derivations}

Let $ \textbf{r}: \Omega \subset \textbf{R}^{2} \longrightarrow \textbf{R}^{3}$ be a parametric representation of a surface $\textbf{S}$: $ \textbf{r}(u,v) $ with its associate normal vector field $\textbf{n}$.

From now on $\textit{E}$, $\textit{F}$, $\textit{G}$ and $\textit{L}$, $\textit{M}$, $\textit{N}$ denote the first and second fundamental forms, respectively, $k_{1}$ and $k_{2}$ the principal curvatures, $H = \frac{k_{1}+k_{2}}{2}$, the mean curvature and $C=\mid k_{1}-k_{2} \mid$ the cylinder.
Given a surface parameterized in orthogonal curvilinear coordinates $(u,v)$, our starting point is Gauss-Codazzi-Mainardi equations for these orthogonal coordinates (p. 149~\cite{Stoker}):
\begin{equation}\label{e1}
\frac{\partial L}{\partial v} =  \frac{\partial E}{\partial v} H,
 \hspace{0.5cm}
\frac{\partial N}{\partial u} =  \frac{\partial G}{\partial u} H.
\end{equation}

Eqs.~(\ref{e1}) are valid except at surface umbilical points, so $k_{1} \neq k_{2}$ must hold.
Now, cylinder can be plugged in Eqs. (\ref{e1}) if one chooses  as orthogonal coordinate system the one associated to curvature lines~\cite{Stoker, Milnor}. Thus:
\begin{equation}\label{cmlc}
C \frac{\partial G}{\partial u} = 2 G \frac{\partial k_{u}}{\partial u}, \hspace{0.5cm}
C \frac{\partial E}{\partial v} = -2 E \frac{\partial k_{v}}{\partial v},
\end{equation}
where the principal curvatures along the lines of curvature $u-v$ are denoted by $k_{1} \equiv k_{v}$ and $k_{2} \equiv k_{u}$.
Note that a regular surface, such as those used in lenses, can be completely covered with lines of curvature. Now, we use the concept of geodesic curvature (the inner curvature of a curve embedded in a surface), which, along a $u-v$ net of curves, is given by (p. 157 \cite{Kreyszig}):
 \begin{equation}\label{kg}
 \begin{split}
&(k_{g})_{u=cte} = \frac{1}{2G\sqrt{E}} \frac{\partial G}{\partial u}, \\
&(k_{g})_{v=cte} = \frac{-1}{2E\sqrt{G}} \frac{\partial E}{\partial v}.
\end{split}
\end{equation}
Substituting Eqs.~(\ref{kg}) into Eqs.~(\ref{cmlc}) we get:
\begin{equation}\label{c}
 C  = \frac{1}{\sqrt{E} (k_{g})_{u=cte}} \frac{\partial k_{u}}{\partial u} = \frac{1}{\sqrt{G} (k_{g})_{v=cte}} \frac{\partial k_{v}}{\partial v}.  \\
\end{equation}
Still, we need to manipulate Eqs.~(\ref{c}) to provide cylinder as a function
of derivatives of curvatures along the curvature lines with the arc-length as curve parameter.
Let us define $s_{u}$ and $s_{v}$ to be the arc length parameters along the $u-v$ net of curves.
For the functions $k_{u}(u,v)$ and $k_{v}(u,v)$, composition derivative rules (and considering that $du/ds_{v}=0$ and $dv/ds_{u}=0$) provide the following equations:
\begin{equation}\label{cd}\begin{split}
&\frac{d k_{u}}{d s_{u}} = \frac{d u}{d s_{u}} \frac{\partial k_{u}}{\partial u}  = \frac{1}{\sqrt{E}} \frac{\partial k_{u}}{\partial u},  \\
&\frac{d k_{u}}{d s_{v}} = \frac{d v}{d s_{v}} \frac{\partial k_{u}}{\partial v}  = \frac{1}{\sqrt{G}} \frac{\partial k_{u}}{\partial v},   \\
&\frac{d k_{v}}{d s_{u}} = \frac{d u}{d s_{u}} \frac{\partial k_{v}}{\partial u}  = \frac{1}{\sqrt{E}} \frac{\partial k_{v}}{\partial u},\\
&\frac{d k_{v}}{d s_{v}} = \frac{d v}{d s_{v}} \frac{\partial k_{v}}{\partial v}  = \frac{1}{\sqrt{G}} \frac{\partial k_{v}}{\partial v}.
\end{split}\end{equation}
Applying Eqs. \eqref{cd} into \eqref{c} provides:
\begin{equation}\label{c2}
 C  = \frac{1}{(k_{g})_{u=cte}} \frac{d k_{u}}{d s_{u}} = \frac{1}{(k_{g})_{v=cte}} \frac{d k_{v}}{d s_{v}}.  \\
\end{equation}

Expression~(\ref{c2}) has been derived somewhat differently elsewhere \cite{Ando, Martin}.

%

Looking in detail to~\eqref{c2} we can extract some relevant conclusions. Cylinder depends linearly on the speed of change of a principal curvature
along one line of curvature and also, at the same time, it is inversely proportional to the geodesic curvature along the orthogonal line of curvature at the same point.
The dependence with curvature properties along both lines of curvatures is an expected result considering that the cylinder is a function of the two principal curvatures, at a specific point of the surface.
Particularly, choosing the principal line to be one of the lines of curvatures we can extract a very important property; contrary to what is predicted by Alonso version of Minkwitz theorem~\cite{Barbero},~\eqref{c2} reveals that there is one, \emph{and only one}, extra \emph{independent} variable for controlling cylinder besides the ratio of change in the principal curvature along the principal line: the geodesic curvature.

This property is not meaningless from a practical point of view, since the geodesic curvature is a measure of how the curve within the surface twists out locally~\cite{Barbero}, which is related to the manufacturing process of that surface. Given a ratio of change of one the principal curvatures, the only way to minimize cylinder is at the expense of \textit{twisting} the surface locally,
which makes manufacturing more complex. In turn, this proves why a tendency to local axial symmetry (along the normal to the surface) is a bad choice for keeping cylinder under control, because in that case the geodesic curvature takes very small values.

We point out that \eqref{c2} becomes undetermined if a line of curvature is coincident with a geodesic line at the evaluation point.
At geodesics the geodesic curvature vanishes, but considering that the cylinder of a smooth function is, itself, smooth, at those points, then the derivative of the principal curvature must also vanish, hence making \eqref{c2} undetermined.
At these points the indeterminacy can be solved applying L'H\^opital's rule along one curvature line  (say, for instance, along $s_u$ if we look  at the first equality from \eqref{c2}). Then we get:
\begin{equation}\label{c2c}
C=\frac{1}{\frac{d}{ds_u}(k_g)_{u=cst}} \frac{d^2 k_u}{ds_u^2}
\end{equation}
at the point under consideration.



Lines of curvature are computed by solving two-first order differential
equations given an initial point $(u_0, v_0)$ on the surface $S$ \cite{Stoker}. Existence and uniqueness of these equations are guaranteed for smooth surfaces (continuity of third derivatives ~\cite{Stoker}), something that is expected in PALs. Nevertheless, we warn that solutions become indeterminate if the lines cross umbilical, planar points (singular points) \cite{Stoker} and points where principal directions are tangent to the surface parameter curves \cite{Farouki}. These points can be avoided by considering that their location is normally known \textit{a priori} in PALs. Indeed, some early designs were based on \textit{a priori} location of a single umbilical line~\cite{Alonso}.

Once the line of curvatures are integrated, they can be reparametrized with respect to the arc-length: $P: [s_i \smallskip s_f] \longrightarrow \Gamma \subset S$. This is easily done interpolating the curve $C$ on a uniform set of points belonging to $[s_i \smallskip s_f]$.
Then, the geodesic curvature can be computed using the formula~\cite{Kreyszig}:
$$k_g  = \left(\textbf{n} \times \frac{dP(s)}{ds}\right) \cdot \frac{d^2 P(s)}{ds^2}.$$
%

Special care must be taken because, whereas the derivative of the principal curvature must be evaluated along one line of curvature, the associated geodesic curvature is numerically obtained evaluating finite differences along the orthogonal line of curvature at the point of interest.

There are available equations for derivatives of normal curvatures along surface parameter curves ~\cite{Mehlum}. Or, alternatively these can be computed numerically by finite differences using the arc-length parametrization of the lines of curvature.
For the numerical evaluation of the cylinder at points where the geodesic curvature vanishes \eqref{c2c} can be evaluated
simply applying finite differences to the derivatives: $\frac{d}{ds_v}(k_g)_{v=cst}$ and $\frac{d^2 k_v}{ds_v^2}$.
%
%
\begin{figure}[h]
\centering
\fbox{\includegraphics[width=\linewidth]{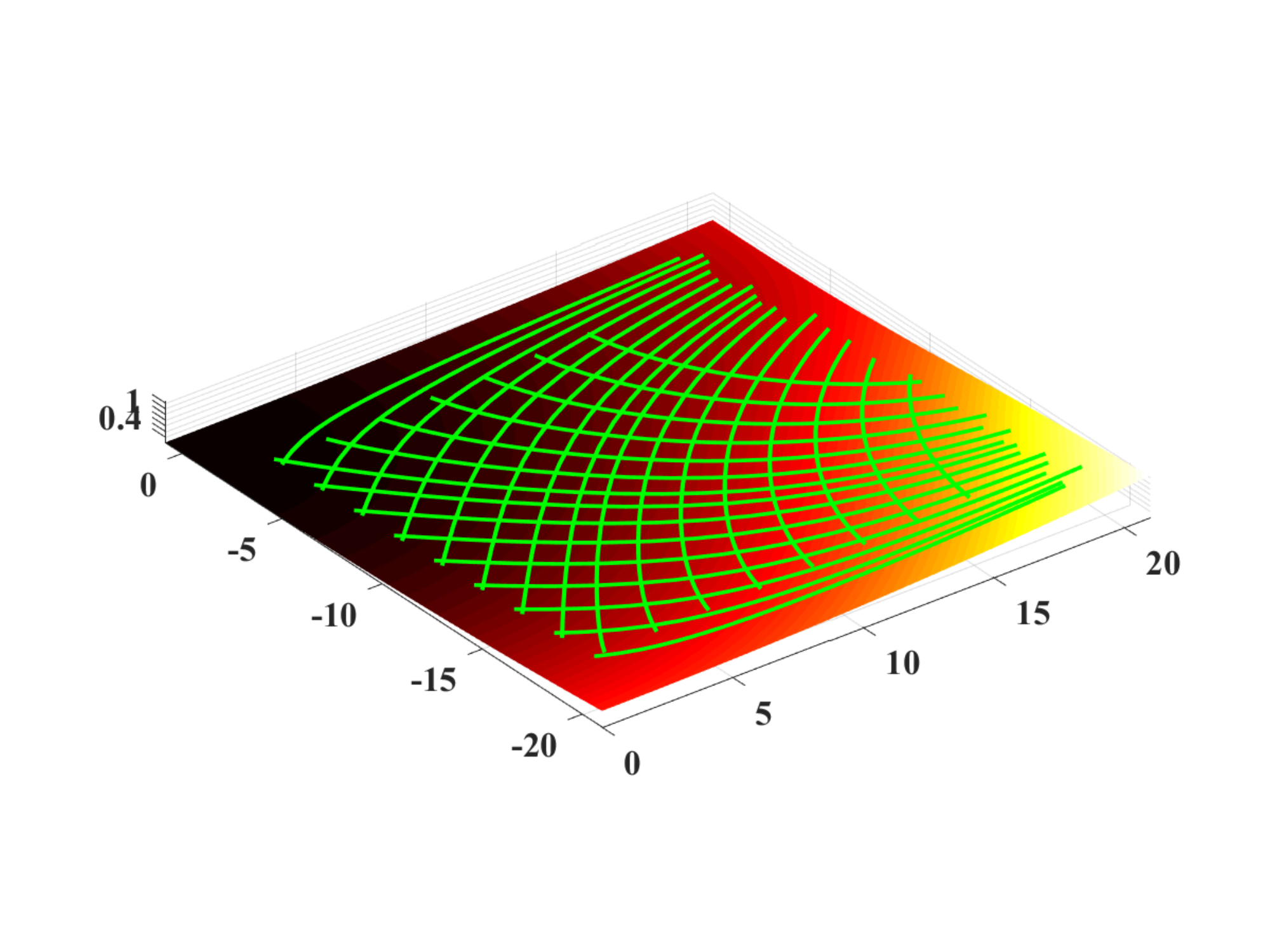}}
\caption{Lines of curvatures net (green lines) associated to the PAL surface given by \eqref{palmodel}. Due to symmetry we only plot the right lower quadrant.}
\label{fig:linescurvature}
\end{figure}

From \eqref{c2} we derive a generalization of Minkwitz theorem to points along lines of curvature and without the low-curvature approximation. Deriving \eqref{c2} with respect to one of the arc-length parameters along the line of curvature and simplifying we obtained:
\begin{equation}\label{gmt}
\begin{split}
 \frac{d C}{d s_{u}}  &= \frac{1}{(k_{g})_{v=cte}} \left(\frac{d^{2} k_{v}}{ds_{u} ds_{v}} -  \frac{1}{(k_{g})_{v=cte}}\frac{d (k_{g})_{v=cte}}{ds_u} \frac{dk_{v}}{ds_v} \right) \\
&=  \frac{1}{(k_{g})_{v=cte}} \left(\frac{d^{2} k_{v}}{ds_{u} ds_{v}} - C \frac{d (k_{g})_{v=cte}}{ds_u} \right).   \\
 \end{split}
\end{equation}
As before, when the geodesic curvature vanishes, \eqref{gmt} could be solved applying L'H\^opital's rule. In contrast to previous Minkwitz generalizations~\cite{Esser, Blendowske}, \eqref{gmt} is exact and valid for any point (except if it is umbilic). It also reveals that the cylinder first derivatives depend on cylinder itself (something already found in~\cite{Esser, Blendowske}) but also on other terms such as geodesic curvature and its first derivative.

\section{Numerical example of application}

To illustrate the implications of the derived compatibility equations, we applied them to an archetypal model of a PAL surface, where the curvature profile along the principal line ($y=0$), under the low-curvature approximation, is cubic-type ~(p. 294~\cite{Alonso}):
\begin{equation}\label{palmodel}
\begin{aligned}
& z(x, y<0)  =  \frac{k_{0}(x^2+y^2)}{2} \\
& -\frac{k_{A}x^2}{2} \left(\frac{10y^3}{L^3} + \frac{15y^4}{L^4} +  \frac{6y^5}{L^5}\right) - k_{A}L^2 \left(\frac{y^5}{2L^5}+\frac{y^6}{2L^6}+\frac{y^7}{7L^7}\right),
\end{aligned}
\end{equation}
where we chose $k_{0} = 2$ D, $k_{A} = 2$ D and $L = 0.002$ m.
Therefore, surface given by ~\eqref{palmodel} provides a variation of 2 D in the mean curvature from far $(y=0)$ to near $(y=-2)$ mm vision. It also contains an umbilical line along the principal line.
\begin{figure}[h]
\centering
\fbox{\includegraphics[width=\linewidth]{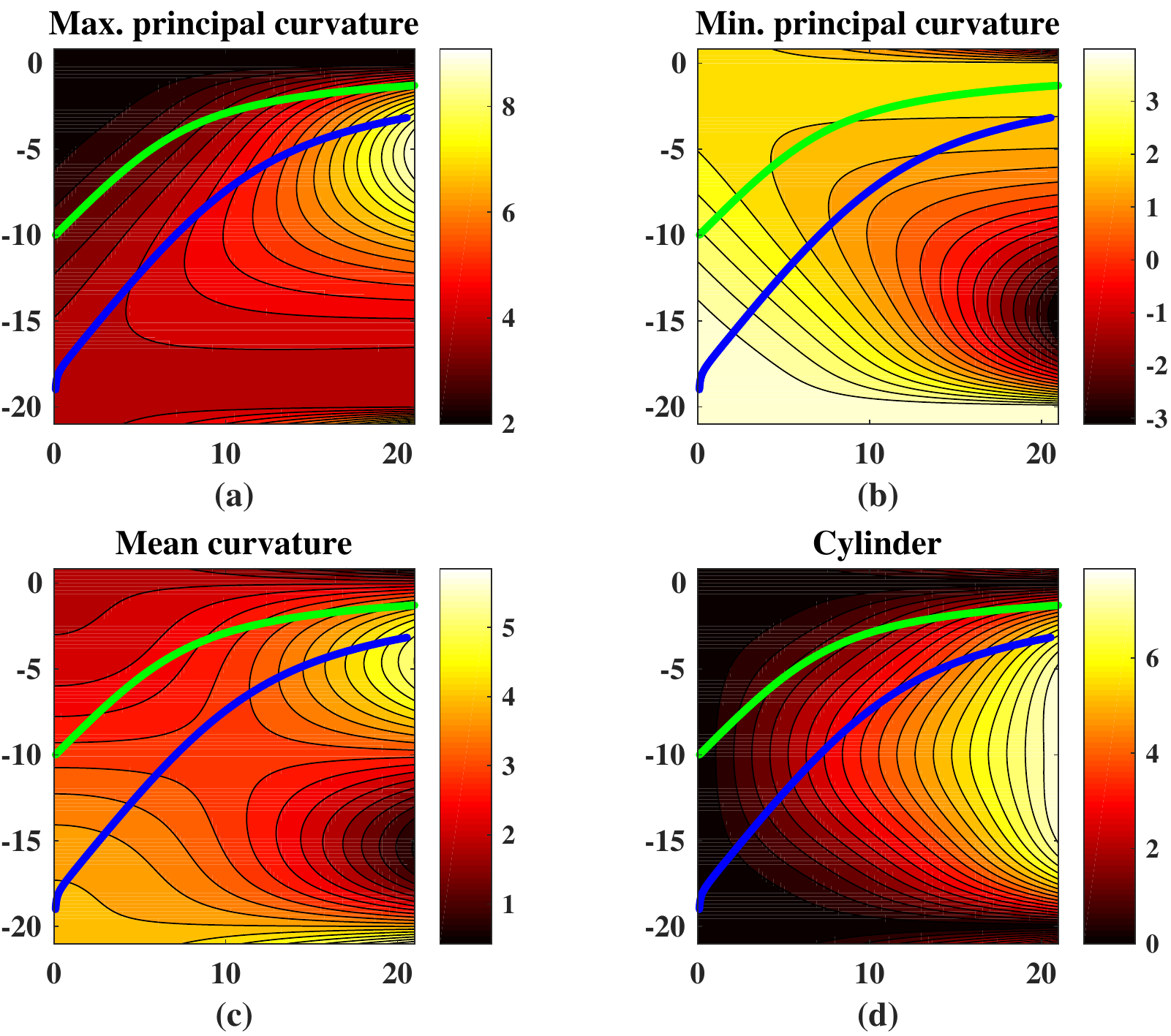}}
\caption{Two lines of curvature (green and blue lines) drawn on: (a) Maximum principal curvature; (b) Minimum principal curvature; (c) Mean curvature; and (d) Cylinder; all of them in diopter units}
\label{fig:linescurvatureonplots}
\end{figure}
We computed the lines of curvatures  within the domain $x \in (0, 20)$ mm
(for $x<0$ the result is symmetrical) and $y \in (-20, 0)$ mm.
The surface and the embedded lines of curvature net (green lines) are shown
in Fig.\ref{fig:linescurvature}.
We note the important fact that the lines of curvature follow special patterns when approaching the umbilical line~\cite{Garcia2005}, sometimes cutting it perpendicular and sometimes not. This implies that not at all the points of the umbilical line ($y=0$) it is possible to associate it, at the limit, to a line of curvature.

\begin{figure}[h]
\centering
\fbox{\includegraphics[width=\linewidth]{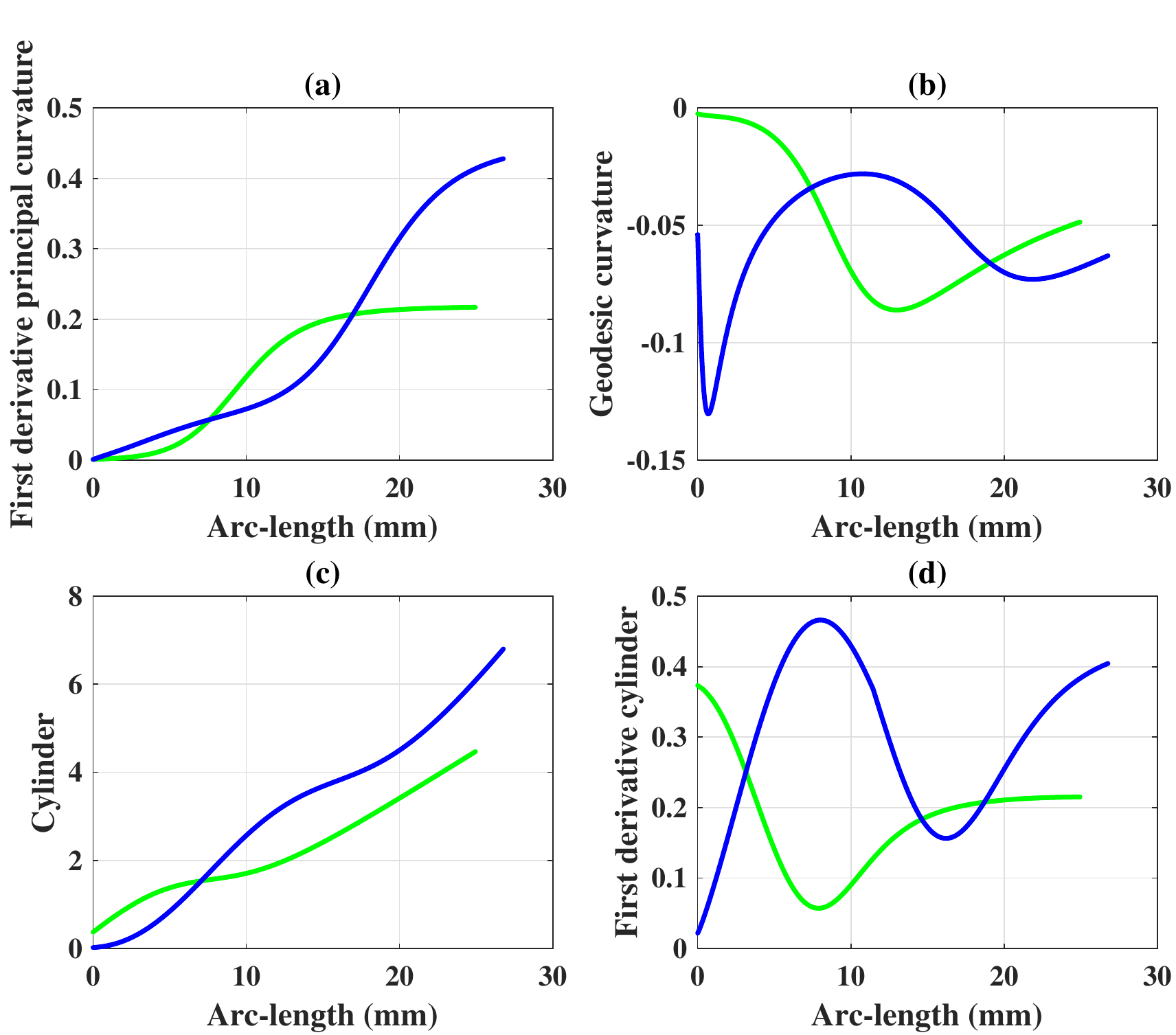}}
\caption{(a) First derivative of maximum principal curvature; (b) Geodesic curvature; (c) Cylinder computed through~\eqref{c2} and the numerical procedure explained in the letter; (d) First derivative of the cylinder with respect to arc-length.
Green and blue curves are lines of curvature for $(1, -10)$ and $(1, -19)$ mm starting points, respectively.}
\label{fig:curves}
\end{figure}

We carried out a more detail analysis in two reference lines of curvature.
First, the line of curvature associated to the maximum principal curvature generated from starting point $(1,-10)$ mm, and second from $(1,-19)$ mm, both increasing along the
positive $x$-axis. Figures \ref{fig:linescurvatureonplots}. show the curves trajectories on the bi-dimensional plots of mean curvature, cylinder and principal curvatures within the selected domain. Both, the principal and the mean curvatures and the cylinder are computed from the principal forms and normals of the surface. As seen in this graph, the two lines cover a wide range of possible values taken by the cylinder.

We analyzed the shape of the first derivative of the maximum principal curvature with respect to arc-length and the geodesic curvatures along the orthogonal curves and its contribution to the cylinder.
Figure \ref{fig:curves} shows first derivative of maximum principal curvature, geodesic curvature and the cylinder and its first derivative with respect to the arc-length as computed through Eqs. (\ref{c2}) and the numerical procedure explained above. The graphs clearly show that the geodesic curvature plays a quite relevant role in the final shape of the form of the cylinder and the cylinder derivative. For instance, in the green line of curvature the region where the geodesic curvature is minimum (around $s=10$ mm)
is close to the maximum of the first derivative of the cylinder.

\begin{figure}[h]
\centering
\fbox{\includegraphics[width=\linewidth]{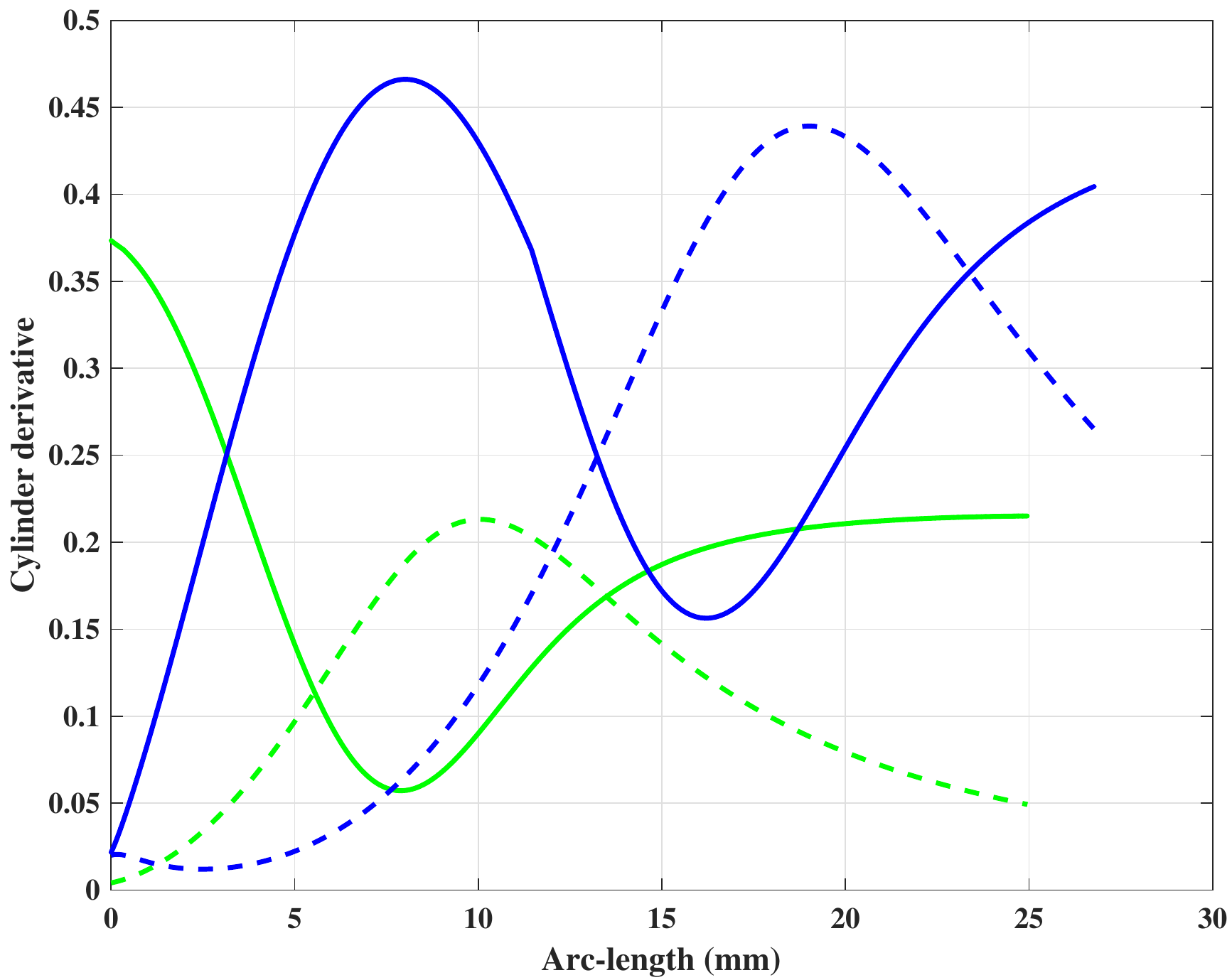}}
\caption{ First derivative of the cylinder with respect to arc-length using
~\eqref{c2} and finite differences on it (solid lines) and Minkwitz theorem (dashed lines).
Green and blue curves are lines of curvature for $(1, -10)$ and $(1, -19)$ mm starting points, respectively.}
\label{fig:ourmink}
\end{figure}

Finally, we compared the cylinder derivative along
the aforementioned lines of curvature numerically obtained by ~\eqref{c2} and finite differences on it with the classical Minkwitz theorem in Fig. \ref{fig:ourmink}.
As expected, the results do not match due to the approximations implied in the derivation of Minkwitz theorem and, also, because Minkwitz is strictly
only valid at umbilical points. However, both equations outputs converge to the same value at the limit, when one of the lines of curvature goes, also at the limit, along the umbilical line, hence matching the $x-y$ axes with the lines of curvature. This is indeed the case of the point obtained intersecting the umbilical line and the line of curvature generated at $(1,-19)$ mm (blue line).

\section{Conclusions}
In summary, we have derived complete and exact compatibility equations, involving cylinder, for a PAL surface to exist.
These equations show that cylinder depends not only on principal curvature derivatives along a line of curvature but also, and exclusively, on the geodesic curvature along an orthogonal line of curvature. This extra degree of freedom is revealed thanks to using lines of curvatures as the surface parametrization.

We have also showed, by numerical examples, that the role
of the geodesic curvature is far from being negligible.
Moreover, we believe that the geometrically
appealing concept of the geodesic curvature, and its possible
relationship with manufacturing processes, could introduce new ways to shed light into understanding the possibilities and restrictions in PAL lens technology.

To end, we would like to note that, of course, optical properties of a PAL surface do not exclusively depend on the surface but also on the geometry of the incoming bundle of light, or wavefront. Following this line, Rubinstein has tried to derived analogous expressions of the restricted Minkwitz theorem to refracted wavefronts~\cite{Rubinstein}. One may try to apply that approach to our new derived expressions.

\section{Funding}
Spanish Ministerio de Economia y Competitividad (FIS2016-75891-P and MTM2017-85757-P).


\bibliographystyle{abbrv}

\end{document}